# An investigation of the asymptotic behavior of the Mertens function

Victor Volfson

ABSTRACT We show in the paper why the asymptotic estimate- the law of the iterated logarithm, does not hold for the Mertens function. A refined asymptotic estimate for the Mertens function is proved.

Keywords: arithmetic function, Mertens function, asymptotic estimate of an arithmetic function, law of the iterated logarithm, probability space, analogue of the law of large numbers, equivalent formulation of the Riemann conjecture.

## 1. INTRODUCTION

The asymptotic behavior of the Mertens function is directly related to the distribution of prime numbers. The Mertens function has the form:

$$M(n)=\sum_{i=1}^{n}\mu(i), \tag{1.1}$$

where $\mu(i)$ is the Möbius function. If a natural has an even number of prime divisors of the first degree, then $\mu(i)=1$, if a natural has an odd number of prime divisors of the first degree, then $\mu(i)=-1$, if a natural has prime divisors of a degree higher than the first, then $\mu(i)=0$.

The following asymptotic estimate for the Mertens function is known:

$$M(n)=O(n^{1/2+\xi}), \tag{1.2}$$

where $\xi$ is a small positive number. The asymptotic estimate (1.2) is an equivalent formulation of the Riemann conjecture [1].

However, there are other asymptotic estimates for the Mertens function. In particular, it was noted in [2, 3] that the Mertens function obeys the law of the iterated logarithm, i.e., almost everywhere the estimate is fulfilled:

$$\lim_{n\to\infty} \sup |M(n)| / \sqrt{n \log\log n} = C, \qquad (1.3)$$

where $C$ is some constant.

The asymptotic estimate (1.3) is true if the Mertens function is a simple symmetric random walk.

A simple random walk is understood as a homogeneous Markov chain $X = (\mathrm{X}_n), n \geq 0$, which describes the movement of a "particle" over integers $Z = \{0, \pm 1, \pm 2, ...\}$, in which this "particle" remains in some state with some probability or can go into one of the neighboring states with some probability. If the transition probabilities to neighboring states are equal, then such a simple random walk is symmetric [4].

Based on formula (1.1), the Mertens function takes only integer values $Z = \{0, \pm 1, \pm 2, ...\}$. At the same time, with some probability, it does not change its value, i.e. remains in the same state or with some probability can go to one of the neighboring states. For the number of steps $\mathrm{n} \to \infty$, the probabilities of transition to neighboring states function are equal for the Mertens.

However, the Möbius function is a multiplicative arithmetic function. Therefore, for any coprime numbers $a, b$, the following holds:

$$\mu(n) = \mu(ab) = \mu(a)\mu(b). \qquad (1.4)$$

Having in mind that the Mertens function can be represented as:

$$\mathrm{M(n)} = \sum_{i=1}^{n-1} \mu(i) + \mu(n), \qquad (1.5)$$

then, based on (1.4) and (1.5), the properties of randomness and independence from prehistory are not satisfied for the Mertens function, i.e. basic properties of a Markov chain.

Thus, the Mertens function is not a simple random symmetric walk. Consequently, the Mertens function does not obey the law of the iterated logarithm, and therefore its asymptotic estimate needs to be refined.

The purpose of this paper is to refine the asymptotic of the Mertens function.

## 2. ASYMPTOTIC ESTIMATE FOR THE ARITHMRTIC FUNCTION

An arithmetic function is a function defined on the set of natural numbers and taking values on the set of complex numbers. The name arithmetic function is due to the fact that this function expresses some arithmetic property of the natural series.

Any initial segment of the natural series $\{1,2,...,n\}$ can be naturally transformed into a probability space $(\Omega_n, \mathcal{A}_n, \mathbb{P}_n)$ by taking $\Omega_n = \{1,2,...,n\}$, $\mathcal{A}_n$— all subsets of $\Omega_n$, $\mathbb{P}_n(A) = \frac{|A|}{n}$. Then an arbitrary (real) function $f(m)$ of the natural argument (more precisely, its restriction to $\Omega_n$) can be considered as a random variable $\xi_n$ on this probability space: $\xi_n(m) = f(m)$, $1 \leqslant m \leqslant n$.

On this probability space, we can talk about the mean value of the arithmetic function $f(m)$:

$$E[f,n] = \frac{1}{n}\sum_{m=1}^{n} f(m), \qquad (2.1)$$

and variance of the arithmetic function $f(m)$:

$$D[f,n] = \frac{1}{n}\sum_{m=1}^{n} |f(m)|^2 - E^2[f,n] \ . \qquad (2.2)$$

Therefore, one can write Chebyshev's inequality for an arithmetic function $f(m), m = 1,...,n$ on this probability space:

$$P_n(|\mathrm{f(m)} - \mathrm{E[f,n]}| \le b\sqrt{D[f,n]}) \ge 1 - 1/\mathrm{b}^2, \qquad (2.3)$$

where $b \ge 1$ .

Let us put $b = b(n)$ in (2.3), where $b(n)$ is an indefinitely increasing function when the value $\mathrm{n} \to \infty$ and get:

$$P_n(|\mathrm{f(m)} - \mathrm{E[f,n]}| \le b(n)\sqrt{D[f,n]}) \ge 1 - 1/\mathrm{b}^2(n). \qquad (2.4)$$

The limit $P_n$ exists in (2.4) when the value $\mathrm{n} \to \infty$:

$$P_n(|\,\mathrm{f(m)} - \mathrm{E[f,n]}\,| \le b(n)\sqrt{D[f,n]}) \to 1, n \to \infty. \quad (2.5)$$

Expression (2.5) is an analogue of the law of large numbers for an arithmetic function [5].

Based on (2.5), we can write an expression for the asymptotic of the arithmetic function $f(m), m = 1,...,n$, which is satisfied almost everywhere when the value $\mathrm{n} \to \infty$:

$$f(n) = E[f,n] + O(b(n)\sqrt{D[f,n]}), \quad (2.6)$$

where $b(n)$ is a slowly growing function, i.e. $\mathrm{b(n)} \le \mathrm{n}^{\xi}$ , a $\xi > 0$.

Having in mind (2.6), if $b(n)\sqrt{D[f,n]} = o(E[f,n])$, then almost everywhere when the value $\mathrm{n} \to \infty$:

$$\mathrm{f(n)} = \mathrm{E[f,n](1+o(1))}, \quad (2.7)$$

otherwise

$$f(n) = O(b(n)\sqrt{D[f,n]}) \ . \quad (2.8)$$

when the value $\mathrm{n} \to \infty$.

Thus, to determine the asymptotic of an arithmetic function $f(m), m = 1,...,n$, it is necessary to know the values $\mathrm{E[f,n]}$ and $D[f,n]$ when the values $\mathrm{n} \to \infty$.

3. DEFINITION OF THE REFINED ASYMPTOTIC OF THE MERTENS FUNCTIONS

QEIS A028442 gives the sequence of zeros of the Mertens function. Let us denote the sequence of these zeros: $z_1, z_2, \ldots$. Segments $[0, z_1], [0, z_2], \ldots$ are nested in this case. The sequence $[0, z_n]$ covers the entire natural series when the value $\mathrm{z}_n \to \infty$.

Based on the above, the equality of probabilities on the segment $[0, z_n]$ is fulfilled:

$$\mathrm{P}(\mu(\mathrm{k}) = 1) = \mathrm{P}(\mu(\mathrm{k}) = -1), \quad (3.1)$$

where $\mathrm{k} \le \mathrm{z}_n$ and $\mathrm{z}_n \to \infty$.

If we consider only $k$, the square-free values, then for the Moebius function truncated in this way, based on equality (3.1), we obtain:

$$P(\mu^*(k)=1)=P(\mu^*(k)=-1)=1/2, \qquad (3.2)$$

where $k \le z_n$ and $z_n \to \infty$.

From (3.2) we obtain that the average value of the modified Mertens function $M^*(n)=\sum_{k=1}^{n}\mu^*(k)$ is equal to:

$$E[M^*,n]=0 \qquad (3.3)$$

when the value $n \to \infty$.

This corresponds to the results [3].

Having in mind (3.3) and that the values where the Möbius function $\mu(k)=0$ does not affect the average value of the Mertens function, we obtain:

$$E[M,n]=0 \qquad (3.4)$$

when the value $n \to \infty$.

Now we define the variance of the Mertens function.

Based on [3] we get:

$$D[M,n]=n \qquad (3.5)$$

when the value $n \to \infty$.

We substitute (3.4) and (3.5) into (2.6) and obtain that, almost everywhere, the following asymptotic estimate holds for the Mertens function:

$$M(n)=O(n^{1/2+\xi}), \qquad (3.6)$$

which corresponds to (1.1) and is an equivalent formulation of the Riemann conjecture.